\documentclass[12pt]{article}
\usepackage{amsmath}
\usepackage{amssymb}
\usepackage{latexsym}
\usepackage{enumerate}

\newtheorem{theorem}{Theorem}[section]
\newtheorem{proposition}[theorem]{Proposition}
\newtheorem{corollary}[theorem]{Corollary}
\newtheorem{lemma}[theorem]{Lemma}
\newtheorem{example}[theorem]{Example}
\newtheorem{remark}[theorem]{Remark}

\newtheorem{definition}[theorem]{Definition}
\newenvironment{proof}{\noindent{\sc Proof.}}{\hfill\qed}

\newcommand{\Z}{{\mathbb Z}}
\newcommand{\Zpinf}{{\mathbb Z_{p^\infty}}}
\newcommand{\F}{{\cal F}}
\newcommand{\calH}{{\cal H}}

\newcommand{\calP}{{\cal P}}

\newcommand{\n}{{\bf n}}

\newcommand{\Hom}{{\rm Hom}}
\newcommand{\qed}{\quad\lower0.05cm\hbox{$\Box$}}

\newcommand{\ra}{\rightarrow}
\newcommand{\lra}{\longrightarrow}

\newcommand{\epi}{\mbox{$\to$\hspace{-0.35cm}$\to$}}

\newcommand{\mono}{\hookrightarrow}
\newcommand{\lemono}{\hookleftarrow}

\def\umono{\ar@{_{(}->}[u]}
\def\uumono{\ar@{_{(}->}[uu]}

\def\lmono{\ar@{_{(}->}[l]}
\def\llmono{\ar@{_{(}->}[ll]}

\begin{document}
\setcounter{page}{1}
\title{Preservation of perfectness and acyclicity:\\
Berrick and Casacuberta's universal acyclic space localized at a
set of primes}

\author{
{\sc Jos\'e L. Rodr\'{\i}guez, J\'er\^ome Scherer, and Antonio
Viruel}
\thanks{Jos\'e L. Rodr\'{\i}guez was partially supported by DGIMCYT grant
BFM2001-2031, and EU grant nr HPRN-CT-1999-0011, J. Scherer by the
program Ram\'on y Cajal, MCYT (Spain), and A. Viruel by MCYT grant
BFM2001-1825, and EU grant nr HPRN-CT-1999-0011.}}

\date{}
\maketitle

\begin{abstract}
In this paper we answer negatively a question posed by
Casacuberta, Farjoun, and Libman about the preservation of perfect
groups under localization functors. Indeed, we show that a certain
$P$-localization of Berrick and Casacuberta's universal acyclic
group is not perfect. We also investigate under which conditions
perfectness is preserved: For instance, we show that if the
localization of a perfect group is finite then it is perfect.

\end{abstract}

\section{Introduction}

The question of whether or not perfect groups are preserved under
arbitrary group localization functors goes back to Farjoun and was
raised simultaneously in \cite{Lib99} and \cite{Cas00}. The close
link between perfect groups and ordinary homology provides the
topological motivation. Indeed if $X$ is a space with perfect
fundamental group $H$, then of course $H_1(X; \Z) = 0$ by the
Hurewicz Theorem. One wonders next what happens with the
localizations of $X$. Is it true that their first integral
homology groups are still trivial? Likewise what can be said about
the localizations of an acyclic space? A first step is to analyze
the effect of localization functors in the category of groups,
since any such functor yields a homotopical localization by
applying the classifying space construction.

A positive result is Theorem~\ref{perfect} which claims that if
$G$ is a perfect group and its localization is known to be finite
then it is automatically perfect. Therefore any example of a
non-perfect localization of a perfect group should be an infinite
group. The reader which is not familiar with the general theory of
localization should not worry since we will mostly deal with
standard $P$-localization, where $P$ is a set of primes. The
functor $L_P$ associates in a universal way to any group $G$ a
morphism $G \ra L_P G$ from $G$ to a uniquely $p$-divisible group
for any prime number $p \not\in P$. We ask now:

\medskip

\noindent {\bf Question}: If $G$ is a perfect group, does it imply
that so is $L_P G$?

\medskip

We have seen that we should look for large groups if we want to
answer negatively to this question. Berrick and Casacuberta have
precisely constructed in~\cite{BC} a universal acyclic space
$B\F=K(\F,1)$, whose associated nullification functor is
equivalent to Quillen's plus construction (cf.\
\cite[1.A.4]{Dro96} for a definition of nullification). The
fundamental group $\F$ is a free product of an uncountable set of
locally free perfect groups and is a universal ``generator'' of
all perfect groups in the sense that the evaluation map
$$
\F\times \Hom(\F,G)\stackrel{{\rm ev}}{\longrightarrow} G
$$
given by ${\rm ev\,}(x,\psi)=\psi(x)$ is surjective for every
perfect group $G$. In particular for any acyclic group $G$ the
cofibre of the evaluation map $\bigvee_{[B\F, BG]} B\F \ra BG$ is
simply connected and acyclic, thus contractible. By Chach\'olski's
result \cite[Theorem~20.9]{Cha} it means that the cellularization
$CW_{B\F} BG$ is $BG$, i.e.\ $BG$ can be constructed as a pointed
homotopy colimit of a diagram whose values are all copies of
$B\F$. It is thus natural to start looking for localizations of
this universal acyclic group $\F$, or of $\calH$, the ``simplest"
of its free summands (the definitions of these groups are given in
Section~4). The main result of this paper is the following
(another example of a non-perfect localization of a perfect group,
based on a quotient of the same group $\calH$, has been obtained
independently by B.~Badzioch and M.~Feshbach in \cite{BF}).

\bigskip

\noindent {\bf Theorem \ref{non-perfect}}
{\it Let $p$ be any prime and $P$ its complementary set of primes.
The $P$-localization $L_P \calH$ of the perfect group $\cal H$ is
not a perfect group. In particular, $L_P\F$ is not perfect
either.}

\bigskip

A straightforward consequence is that the topological
$P$-localization of the acyclic space $B\F$ is not acyclic, as its
fundamental group is precisely $L_P \F$.

\bigskip

\noindent {\bf Corollary \ref{classifying}}
{\it Let $p$ be any prime and $P$ its complementary set of primes.
The $P$-localization $L_P B\calH$ of the acyclic space $B \calH$
is not acyclic. In particular, $L_P B\F$ is not acyclic either.}

\bigskip

This result on the first homology group contrasts with the main
theorem in \cite{CS}, which deals with the first homotopy group:
Any homological localization of a simply connected space is again
simply connected.

The present work fits into an already large series of recent
papers aimed at the study of preservation of algebraic structures
by localization functors, with applications to stable and unstable
homotopy theory. Let us mention the nice survey \cite{Cas00}, as
well as \cite{Lib99}, \cite{GRS}, \cite{S}, \cite{RST},
\cite{RSV}, and \cite{CG}.

\medskip

{\it Acknowledgments:} We warmly thank the referee for the
tremendous improvement in the organization of the material of
Section 3 and 4, and in particular for his proof of
Proposition~\ref{nicegroup}. We would also like to thank Warren
Dicks for pointing out Example~\ref{almost} and Carles Casacuberta
for helpful comments.

\section{Finite localizations of perfect groups}

The main result of this section could be useful to determine all
finite localizations of finite simple groups, a problem which was
extensively studied in \cite{RST}. It was shown indeed
in~\cite{RSV} that the localization of a finite simple group may
(un)fortunately not be simple. We prove nevertheless in this
section, Theorem~\ref{simple-perfect}, that all finite
localizations of a non-abelian finite simple group must be perfect
groups.

Let $(L,\eta)$ be any localization functor in the category of
groups, where $\eta: G \ra LG$ denotes the natural coaugmentation.
Recall briefly that a group $G$ is called $L$-local if $\eta: G
\cong LG$, and it is called $L$-acyclic if $LG=1$. Note that for
every $L$-acyclic group $H$ and every $L$-local group $G$, we have
$\Hom(H,G)=0$. We will say that a morphism $\varphi: H \ra G$ is a
localization if there exists a localization functor $L$ such that
$\varphi$ coincides with the coaugmentation $\eta: H \ra LH$. This
actually happens exactly when $\varphi$ induces a bijection
$\Hom(H, G) \cong \Hom(G, G)$, see \cite[Lemma~2.1]{Cas00}.

\begin{definition}
{\rm A group $G$ is called {\it almost perfect} if it has no
non-trivial endomorphism which factors through an abelian group,
or equivalently through $G_{ab}$.}
\end{definition}

We do not know if this property has a name in the literature, but
it makes sense for any variety ${\cal V}$ of groups (for more
about varieties, see for example \cite[Section 1.1]{Rob}). A group
could be called {\it almost ${\cal V}$-perfect}, if it has no
non-trivial endomorphism which factors through a group in ${\cal
V}$. Of course any perfect group is almost perfect. On the other
hand there are non-perfect groups which are almost perfect. Our
main result Theorem~\ref{non-perfect} actually shows that the
group $L_P \calH$ is not perfect but almost perfect. Let us
indicate now a more natural and elementary example of this type.

\begin{example}
\label{almost}
{\rm In \cite{St} Strebel considers the following central
extensions of the triangle groups: $G(l,m,n) = \langle x, y \, |\,
x^l = y^m = (xy)^n \rangle$. They all occur as fundamental groups
of certain 3-dimensional manifolds and are finite if and only if
$|l|^{-1}+|m|^{-1}+|n|^{-1} > 1$. The infinite ones are Poincar\'e
duality groups of dimension 3 (therefore torsion-free) and the
abelianization can be easily computed. For example $G=G(3, 8, 2)$
is such a torsion-free group with $G_{ab} \cong \Z/2$. It is thus
almost perfect, but not perfect.}
\end{example}

\begin{lemma}
Let $H$ be a perfect group and $H \ra G$ be a localization. Then
$G$ is almost perfect.
\end{lemma}

\begin{proof}
Consider an endomorphism $\psi: G\to G$ which factors through an
abelian group~$A$. Since $H$ is perfect, the composite $\psi \circ
\eta: H\to G\to A\to G$ is trivial. On the other hand, there is a
bijection of sets $\eta^*: \Hom(G,G)\cong \Hom(H,G)$ by the
universal property of the localization. As $\eta^*(\psi)=0$, we
infer that $\psi$ is trivial.
\end{proof}

\begin{theorem}
\label{perfect}
Let $H$ be a perfect group and $H \ra G$ be any localization with
$G$ finite. Then $G$ is perfect.
\end{theorem}

\begin{proof}
Assume $G$ is not perfect. Then there exists a prime $p$ such that
$\Z/p$ is a quotient of $G_{ab}$ because $G$ is finite. But then
we would have a non-trivial endomorphism $G \epi G_{ab} \epi \Z/p
\mono G$, which contradicts the assumption that $G$ is almost
perfect.
\end{proof}

\medskip

The following corollary tells us that the kind of constructions we
made in the previous paper \cite{RSV} will never yield an example
of non-perfect localization of a perfect group. We have to deal
with infinite groups.

\begin{corollary}
\label{simple-perfect}
Let $H$ be a non-abelian finite simple group and $H \ra G$ be any
localization of finite groups. Then $G$ is perfect. \hfill{\qed}
\end{corollary}

\begin{remark}
\label{large}
{\rm The work of R.~G\"{o}bel and S.~Shelah \cite{GS} (see also
\cite{GRS}) shows that every finite simple group admits large
localizations (as large as any given cardinal). We do not know if
a similar result also holds for localizations of finite perfect
groups.}
\end{remark}

\section{Baumslag's class $\calP$}

Given a set of primes $P$, a group $G$ is called $P$-local if the
$p$-power map $x\to x^p$ is bijective for all $p\not\in P$. In the
notation of \cite{Bau}, such a group belongs to the class
$D_\omega$. Every group $G$ admits a functorial $P$-localization
$\eta_G: G\to L_PG$, i.e.\ $L_PG$ is $P$-local and for every
homomorphism $\psi:G\to K$ into a $P$-local group $K$ there exists
a unique homomorphism $\tilde \psi: L_PG\to K$ such that $\tilde
\psi \circ \eta_G=\psi$. From now on we will invert a single
prime~$p$. Thus $P$ is the set of all primes different from $p$
and the localization $G \ra L_P G$ is obtained by adding unique
$p$-roots to elements without $p$-roots and identifying $p$-roots
when they are not unique. In general $L_P G$ is a very large
group, where the equation $y=x^p$ has a unique solution.

We shall suppose that the reader is somewhat familiar with
Baumslag's work \cite{Bau} on groups with unique roots. Another
valuable reference is Ribenboim's article \cite{Rib}. An important
class of groups introduced by Baumslag in \cite[Section 27]{Bau}
is the class $\calP$ ($\calP_\omega$ in Baumslag's notation). For
a group $G$ in this class he shows that
$$
L_P G = \bigcup_{\alpha < \lambda} G_\alpha
$$
where $\lambda$ is a limit ordinal and $\{G_\alpha \}_\alpha$ is a
chain of supergroups of $G$, all of them belonging to the class
$\calP$. Moreover, for every ordinal $\alpha < \lambda$, either
$G_{\alpha^+} = G_\alpha$ or $G_{\alpha^+}$ is constructed as the
push-out of a diagram
$$
\Z[1/p] \lemono \Z \stackrel{x_\alpha}{\lra} G_\alpha,
$$
i.e.\ by adding $p$-roots to an element $x_\alpha$ in $G_\alpha$.
If $\alpha$ is a limit ordinal then $G_\alpha = \cup_{\beta <
\alpha} G_\beta$. This construction is explained by Baumslag in
\cite[(33.3)]{Bau}. In particular groups in the class $\calP$ are
contained in their $P$-localization: $G \subseteq L_P G$. The
importance of the class $\calP$ comes from
\cite[Corollary~35.7]{Bau}, which shows that free groups belong to
this class. We proceed by proving a lemma, the proof of which
resembles the argument of \cite[Theorem~33.4]{Bau}.

\begin{lemma}
\label{PnotD}
Let $G$ be a group in $\calP$ which is not $P$-local. Then there
exists an epimorphism $L_P G \epi \Zpinf$.
\end{lemma}

\begin{proof}
Since $G$ is not $P$-local there is a minimal $\alpha$ such that
$G_\alpha \neq G_{\alpha^+}$. Therefore $G_{\alpha^+}$ is the
colimit of a certain push-out diagram
$$
\Z[1/p] \lemono \Z \stackrel{x_\alpha}{\lra} G_\alpha.
$$
Construct now an epimorphism $G_{\alpha^+} \epi \Zpinf$ by mapping
$G_\alpha$ trivially to $\Zpinf$ and $\Z[1/p]$ to $\Zpinf$ by the
canonical projection.

The proof goes on by a transfinite induction. Assume the
epimorphism $G_\beta \epi \Zpinf$ has already been constructed for
some ordinal $\beta > \alpha$. We need to check that it can be
extended to the next step $G_{\beta^+}$. If $G_\beta =
G_{\beta^+}$ there is nothing to do, so we may assume that
$$
G_{\beta^+} = colim(\Z[1/p] \hookleftarrow \Z \ra G_\beta).
$$
A morphism out of it is determined by a pair of compatible
morphisms out of $\Z[1/p]$ and~$G_\beta$. Consider the composite
$\Z \ra G_\beta \ra \Zpinf$. The image of $1 \in \Z$ generates a
cyclic subgroup of order $p^k$ in the Pr\"ufer group $\Zpinf$, for
some $k \geq 0$. We require now the morphism $\Z[1/p] \ra \Zpinf$
to be the projection $\Z[1/p] \epi \Z[1/p]/p^k\Z \cong \Zpinf$.
This defines the desired epimorphism $G_{\beta^+} \epi \Zpinf$ and
the inductive process yields then an epimorphism $L_P(G) \epi
\Zpinf$.
\end{proof}

\bigskip

\section{A counterexample}

Berrick and Casacuberta constructed in \cite{BC} an acyclic
two-dimensional space $B\F=K(\F,1)$ whose associated nullification
functor is equivalent to Quillen's plus construction. The ideas
behind its construction have been used in other contexts, for
instance to show the existence of a plus-construction in the
category of differential graded algebras over an operad
\cite{ChRS}. The key ingredient is the universal acyclic group
$\F$, which is a free product of an uncountable set of locally
free perfect groups $\F_{\n}$ indexed by non-decreasing sequences
$\n=(n_1, n_2, \cdots)$ of positive integers. We will make use of
the group $\calH = \F_{(1,1,1, \cdots)}$. The presentation of
$\calH$ is by definition
$$
\langle\, x_1, x_2, x_3, \, \dots \, | \, x_1=[x_2,x_3], \,
x_2=[x_4,x_5], \, x_3=[x_6,x_7], \, x_4 = [x_8, x_9], \, \dots \,
\rangle,
$$
so every generator is a simple commutator of two new generators.
In other words, $\calH$ is the colimit of the diagram
$$
F_1 \stackrel{\varphi_0}{\longrightarrow} F_2
\stackrel{\varphi_1}{\longrightarrow} F_4
\stackrel{\varphi_2}{\longrightarrow} F_8
\stackrel{\varphi_3}{\longrightarrow} \cdots
$$
where $F_{2^n}$ is a free group on generators $x_{2^n}, x_{2^n
+1}, \dots, x_{2^{n+1}-1}$ and $\varphi_{n}(x_i) = [x_{2i} ,
x_{2i+1}]$ for any $2^n \leq i < 2^{n+1}$ (thus $\varphi_{n}$ is
actually a free product of $2^n$ copies of $\varphi_0$). In
general the group $\F_{(n_1,n_2, \cdots)}$ is defined similarly:
the generators at the step $k$ are products of $n_k$ commutators
of new generators; for more details see \cite[Example~5.3]{BC}.
The rest of the section is devoted to proving that the
$P$-localization of the group $\calH$ is not perfect.

We shall need the following terminology from \cite{CGT}. Elements
of the free group $F_n$ on $n$ generators $x_1, \, \dots , x_n$
correspond to freely reduced words $w$ in the symbols $x_1, \,
\dots , x_n$: In $w$ the symbols $x_i^{\epsilon}$ and
$x_i^{-\epsilon}$ never appear adjacently, where $\epsilon = \pm
1$. A word is called cyclically reduced if it is not of the form
$rur^{-1}$ for some words $r$ and $u$.

\begin{lemma}
\label{words}
Let $F_n$ be the free group on $n$ generators $x_1, \, \dots ,
x_n$ and let $w$ be a freely reduced word. If $w$ belongs to $F_n
\setminus F_{n-1}$, then so does $w^k$ for any $k \neq 0$.
\end{lemma}

\begin{proof}
The proof of \cite[Corollary~1.2.2]{CGT} can be used nearly
verbatim. We assume $k>0$ (the case $k<0$ will follow directly).
If $w$ is cyclically reduced and contains the symbol $x_n$, then
$w^k$ is freely reduced and contains thus the symbol $x_n$ as well
for any $k\neq 0$. Moreover the freely reduced word representing
$w^k$ starts with the same letter and ends with the same letter as
$w$.

If $w$ is not cyclically reduced, it is of the form $w=rur^{-1}$
for some cyclically reduced and non-trivial word $u$. Then $w^k =
ru^k r^{-1}$, which is freely reduced since $u^k$ is so and starts
and ends with the same symbol as $u$. Therefore $w^k$ contains the
letter $x_n$ as either $u^k$ or $r$ contains it.
\end{proof}

\begin{proposition}
\label{rootimage}
Let $\varphi:F_n \ra F_{2n}$ be the group homomorphism defined by
$\varphi(x_i) = [x_{2i-1}, x_{2i}]$ for any $1 \leq i \leq n$.
Assume $w \in F_n$ is a word such that $\varphi(w) = u^k$ for some
$k\neq 0$ and some element $u \in F_{2n}$. There exists then an
element $v \in F_n$ such that $\varphi(v) = u$ and $v^k = w$. In
particular, if $w$ has no $p$-root, then neither has $\varphi(w)$.
\end{proposition}

\begin{proof}
In \cite[Lemma~1.1]{HW}, Hurley and Ward give an explicit basis
for the commutator subgroup of a free group. In particular
$[F_{2n}, F_{2n}]$ is a free group on a set of elements containing
the commutators $[x_{2i-1}, x_{2i}]$. This shows first that
$\varphi$ is a monomorphism and second that $[F_{2n}, F_{2n}]$
decomposes as $\varphi(F_n) * E$ for some free subgroup $E$ of
$F_{2n}$.

Notice next that $u$ must belong to $[F_{2n}, F_{2n}]$ since $u^k$
does so and $F_{2n}/ [F_{2n}, F_{2n}] = (F_{2n})_{ab}$ is torsion
free. Thus $u$ actually belongs to the image of $\varphi$ by
Lemma~\ref{words}. There exists therefore an element $v \in F_n$
such that $\varphi(v) = u$. As $\varphi$ is injective, $v^k = w$.
\end{proof}

\begin{proposition}
\label{notlocal}
The group $\calH$ is not $P$-local.
\end{proposition}

\begin{proof}
Consider the inclusion $F_1 \mono \calH$. We prove that the image
of $x_1$ in $\calH$ has no $p$-root. It suffices to do so for the
image of $x_1$ in $F_{2^n}$ for all $n \geq 1$. The proof goes by
induction. Clearly $x_1$ has no $p$-root in $F_1$. The induction
step is provided by Proposition~\ref{rootimage}.
\end{proof}

\begin{proposition}
\label{nicegroup}
The group $\calH$ is in the class $\calP$.
\end{proposition}

\begin{proof}
We want to employ \cite[Lemma~31.1]{Bau} to the expression of
$\calH$ as the colimit of
$$
F_1 \stackrel{\varphi_0}{\longrightarrow} F_2
\stackrel{\varphi_1}{\longrightarrow} F_4
\stackrel{\varphi_2}{\longrightarrow} F_8
\stackrel{\varphi_3}{\longrightarrow} \cdots
$$
By \cite[Corollary~35.7]{Bau} each free group $F_{2^n}$ is in the
class $\calP$. It suffices thus to prove that at each stage in
this telescope the image $\varphi_{n}(C_{F_{2^{n}}}(w))$ of the
centralizer of any element $w \in F_{2^n}$ is equal to the
(infinite cyclic) centralizer $C:= C_{F_{2^{n+1}}}(\varphi_n(w))$
of its image.

Let $u$ be a generator of $C$. Since $\varphi_n(w)$ belongs to $C$
and $\varphi_n$ is injective, there exists some $k \neq 0$ such
that $\varphi_n(w) = u^k$. Proposition~\ref{rootimage} ensures now
the existence of an element $v \in F_{2^n}$ with $\varphi_n(v) =
u$ and $v^k = w$. As $\varphi_n$ is injective, the elements $v$
and $w$ already commute in $F_{2^n}$. It follows that $v \in
C_{F_{2^{n}}}(w)$ and so
$$
C_{F_{2^{n+1}}}(\varphi_n(w)) \subset
\varphi_{n}(C_{F_{2^{n}}}(w)).
$$
The opposite inclusion is obvious, and the result follows.
\end{proof}

\begin{theorem}
\label{non-perfect}
Let $p$ be any prime and $P$ its complementary set of primes. The
$P$-localization $L_P \calH$ of the perfect group $\cal H$ is not
a perfect group. In particular, $L_P\F$ is not perfect either.
\end{theorem}

\begin{proof}
We have seen in Proposition~\ref{notlocal} that the group $\calH$
is not $P$-local, but belongs to the class $\calP$
(Proposition~\ref{nicegroup}). We conclude then by
Lemma~\ref{PnotD} that there exists an epimorphism $L_P \calH \epi
\Zpinf$. Thus $L_P \calH$ is not a perfect group.

The group $\F$ contains by definition $\calH$ as a free summand.
Therefore $L_P \calH$ is a retract of $L_P \F$, which cannot be
perfect.
\end{proof}

\medskip

Recall that in the homotopical context $L_P$ denotes
localization with respect to the $p$-power map $p: S^1 \ra S^1$.
The article \cite{CP} contains valuable information about this
$P$-localization.

\begin{corollary}
\label{classifying}
Let $p$ be any prime and $P$ its complementary set of primes. The
$P$-localization $L_P B\calH$ of the acyclic space $B \calH$ is
not acyclic. In particular, $L_P B\F$ is not acyclic either.
\end{corollary}

\begin{proof}
By \cite[Theorem~8.7]{CP} the topological $P$-localization of a
wedge of circles is given by the group theoretical one. That is,
for any free group $F$, we have $L_P BF \simeq  B(L_P F)$. Since
$B\calH$ is constructed as a telescope of $BF_{2^n}$'s,
Theorem~1.D.3 in \cite{Dro96} tells us that $L_P B\calH \simeq L_P
(hocolim L_P BF_{2^n})$. Therefore
$$
L_P B\calH \simeq L_P (hocolim B L_P F_{2^n}) \simeq L_P B(colim
L_P F_{2^n}) \simeq L_P B(L_P \calH).
$$
The last weak equivalence is due to the fact that $L_P \calH$ is
isomorphic to $colim_{n} L_P F_{2^n}$, this colimit being a
$P$-local group (uniquely $p$-divisible) by a standard finite
object argument. We conclude now that $L_P B\calH \simeq B(L_P
\calH)$ because the latter is a $P$-local space (clearly $BG$ is
$P$-local for any $P$-local group $G$).
\end{proof}

\medskip

It could be nice to completely determine the abelianization of the
$P$-localization of $\F$ (we just exhibited one morphism onto a
Pr\"ufer group), and also to compute explicitly the topological
$P$-localization of Berrick and Casacuberta's universal space
$B\F$.


\vskip 0.5 cm

\setlength{\baselineskip}{0.6cm}

\bigbreak

Jos\'e L. Rodr\'{\i}guez

\noindent \'Area de Geometr\'\i a y Topolog\'\i a, CITE III,
Universidad de Almer\'\i a, E--04120 Almer\'\i a, Spain, e-mail:
{\tt jlrodri@ual.es}

\bigskip

J\'er\^ome Scherer

\noindent Departament de Matem\`atiques, Universitat Aut\`onoma de
Barcelona, E--08193 Bellaterra, Spain, e-mail: {\tt
jscherer@mat.uab.es}

\bigskip

Antonio Viruel

\noindent Departamento de Algebra, Geometr\'\i a y Topolog\'\i a,
Universidad de M\'alaga, Apartado de Correos 59, E--29080
M\'alaga, Spain, e-mail: {\tt viruel@agt.cie.uma.es}

\end{document}